\documentclass[12pt]{amsart}
\usepackage{fullpage}
\usepackage[all]{xy}
\numberwithin{equation}{section}

\theoremstyle{definition}

\title{A counterexample for lightning flash modules over $E(e_1,e_2)$}
\author{David Benson and Robert R. Bruner}
\thanks{
This work  was partially supported by the Simons Foundation and 
by the Mathematisches Forschungsinstitut Oberwolfach.}

\begin{document}

\begin{abstract} 
We give a counterexample to Theorem 5 in \S18.2 of Margolis'  
book, ``Spectra and the Steenrod Algebra,'' and make remarks about the proofs  
of some later theorems in the book that depend on it.  
The counterexample is a module which does not split as a sum of lightning 
flash modules and free modules.
\end{abstract}
 
\maketitle
\section{Introduction}

Let $k$ be a field and $E(e_1,e_2)$ be a graded exterior algebra on
generators $e_1$ and $e_2$ with degrees satisfying $0<|e_1|<|e_2|$.
Theorem 5 in \S18.2 of Margolis \cite{Margolis:1983a} states
that every graded $E(e_1,e_2)$-module is a coproduct of free modules 
and lightning flashes. In this note, we give a simple counterexample 
to this statement.

Statement (c) following Proposition 7 of the same
section is true, but not because of Theorem 5. The proof of Theorem 8
in \S18.3 depends on this statement. The proofs of Proposition 9 and
Lemma 10 of the same section also depend on Theorem 5, and are used
in Chapter 20. Fortunately, the paper of Adams and Margolis \cite{Adams/Margolis:1971a}
provides correct proofs of these statements that do not rely on Theorem 5.

\section{The counterexample}

In this section we display a bounded below module $M$ for $E(e_1,e_2)$ 
which is not isomorphic to a coproduct of free modules and lightning flashes.

First we note that every module for $E(e_1,e_2)$ can be written as a direct sum 
of a free module and a module on which $e_1e_2$ acts as zero. So we may
as well work with modules for $E(e_1,e_2)/(e_1e_2)$.

We use the notation of \S18.2 of Margolis. Let $M(n)$ be the lightning module
$L(n,0,1)$ of dimension $2n$. Here is a picture of $M(n)$:
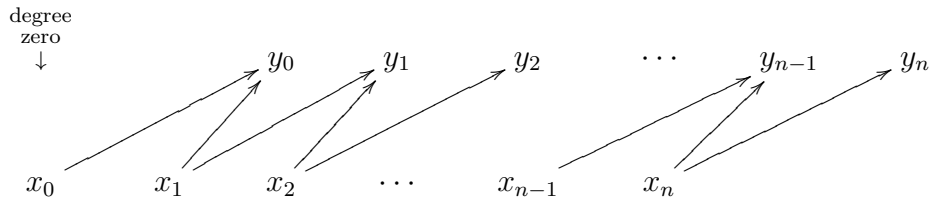
\[ \xymatrix{\shortstack{\tiny degree \\ \tiny zero \\ $\scriptstyle\downarrow$}&&y_0&y_1&y_2&\cdots&y_{n-1}&y_{n} \\ 
x_0 \ar[urr]&x_1 \ar[ur]\ar[urr]&x_2\ar[ur]\ar[urr]&\cdots&x_{n-1}\ar[urr]&x_{n}\ar[urr] \ar[ur]} \]
The shorter arrows represent the action of $e_1$, and the longer ones $e_2$. Thus a presentation of 
the module is given by $e_1x_{i+1}=e_2x_i=y_i$ ($0\le i \le n-1$), $e_1x_0=0$,
$e_2 x_n = y_n$.
We arrange that the element $x_0$ in $M(n)$ is in degree zero, so that $x_i$ has
degree $i(|e_2|-|e_1|)$ and $y_i$ has degree $|x_i|+|e_2|$. 
Similarly, 
$L(\infty,0)$ is the infinite lightning flash obtained by letting this diagram continue
to the right indefinitely.

Our counterexample is the module
\[ M = \prod_{n=0}^\infty M(n). \]
To see that it is a counterexample, first note that $e_1M(n)$ is the linear span of 
$y_0,\dots,y_{n-1}$, so $e_2^{-1}e_1M(n)$ is the linear span of all the basis elements
except $x_n$. Here, if $U$ is a linear subspace of a module, we write $e_2^{-1}U$
for the linear subspace consisting of the vectors whose image under $e_2$ is in $U$. 

Inductively, we see that for $j>0$, $(e_2^{-1}e_1)^jM(n)$ is the linear span of
the basis elements $y_0,\dots,y_n,x_0,\dots,x_{n-j}$. Thus $x_0$ is in $(e_2^{-1}e_1)^jM(n)$
if and only if $j\le n$. 

Taking degree zero parts, we have
\[ ((e_2^{-1}e_1)^jM)_0=\prod_{j=n}^\infty M(n)_0. \]
Thus 
\begin{equation}\label{eq:intersection} 
\bigcap_{j\ge 0}((e_2^{-1}e_1)^jM)_0=0, 
\end{equation}
and
\[ ((e_2^{-1}e_1)^jM)_0/(e_2^{-1}e_1)^{j-1}M)_0\]
is one dimensional.
On the other hand, $x_0$ is in $(e_2^{-1}e_1)^jL(\infty,0)$ for all $j> 0$, so
\[ \bigcap_{j\ge 0}((e_2^{-1}e_1)^jL(\infty,0))_0\ne 0. \]
Since a finite sum is always a direct summand of the product, 
it follows that $M$ has exactly one copy of each $M(n)$ as a summand, and
no summand isomorphic to $L(\infty,0)$. Since $e_1M_0=0$, no summand of
the form $L(\infty,1)$, $L(n,1,0)$ or $L(n,1,1)$ can contribute to $M_0$;
and finally \eqref{eq:intersection} shows that no summand of the form
$L(n,0,0)$ can contribute to $M_0$, since that intersection is non-zero for
such a module.
The summands we have identified do not exhaust $M_0$, and hence
$M$ cannot be a direct  sum of lightning flash  modules.

On the other hand, modules of finite type for $E(e_1,e_2)$ can be shown
to be direct sums of lightning flashes, by the method of filtrations of the
forgetful functor to graded vector spaces. The proof is similar to but easier
than the functorial filtration proof given in Ringel \cite{Ringel:1975a}.

\bibliographystyle{amsplain}
\bibliography{../repcoh}

\end{document}